\documentclass[12pt,psamsfonts]{article}

\usepackage{theorem}
\usepackage{amssymb}
\usepackage{amsmath}

\newtheorem{thm}{Theorem}[section]
\newtheorem{prop}[thm]{Proposition}
\newtheorem{cor}[thm]{Corollary}
\newtheorem{lm}[thm]{Lemma}

{\theorembodyfont{\rm}

\newtheorem{defn}[thm]{Definition}

\newtheorem{rem}[thm]{Remark}

}

\def\ra{\rightarrow}

\def\C{{\Bbb C}}
\def\F{{\Bbb F}}

\def\P{{\Bbb P}}

\def\T{{\bf T}}

\def\F{{\overline{F}}}

\title{On the density of rational points on 
elliptic fibrations}

\author{F. A.  Bogomolov\\
\small  Courant Institute of Mathematical Sciences, N.Y.U. \\
\small 251 Mercer str. \\
\small New York, (NY) 10012, U.S.A.\\
\small e-mail: bogomolo@cims.nyu.edu\\
\small  and\\
Yu. Tschinkel\\
\small Dept. of Mathematics, U.I.C.\\
\small 851 South Morgan str.\\
\small Chicago, (IL) 60607-7045,  U.S.A.  \\
\small e-mail: yuri@math.uic.edu
}

\begin{document}

\date{}

\maketitle

\thispagestyle{empty}

\pagebreak

\section{Introduction}

Let $X$ be an algebraic variety defined over a number field $F$.
We will say that rational points are {\em potentially dense}
if there exists a finite extension $K/F$ such that 
the set of $K$-rational points $X(K)$ is Zariski dense in $X$. 
The main problem is to relate this property to geometric
invariants of $X$.
Hypothetically, on varieties of 
general type rational points are not
potentially dense. In this paper we are interested in 
smooth projective varieties such that neither they
nor their unramified coverings admit
a dominant map onto varieties of general type. 
For these varieties it seems plausible to expect 
that rational points are potentially dense 
(see \cite{harris-tschi-98}). 

\

Varieties which are not of general type can be thought of as 
triple fibrations  $X\ra Y\ra Z$, where  
the generic fiber of $ X\ra Y$ is rationally connected, 
$Y\ra Z$ is a Kodaira fibration with generic fiber 
of Kodaira dimension $0$ and the base has 
Kodaira dimension $\le 0$ (cf. \cite{kollar-miyaoka-mori-92-2}). 
In this paper we study mostly varieties of dimension $2$ and $3$. 

\

In dimension $2$ the picture is as follows:
Rationally connected surfaces are 
rational over some finite extension of $F$ and therefore
the problem has an easy solution in this case. 
Rational points on abelian surfaces are potentially dense. 
In \cite{harris-tschi-98} density was proved for
smooth quartic surfaces in $\P^3$ which contain a line and 
for elliptic fibrations over $\P^1$, provided they have
{\em irreducible} fibers and
a rational or elliptic multisection. 
We don't know the answer for general K3 surfaces.
In this paper we prove that rational points are potentially dense on 
double covers of $\P^2$ ramified in a {\em singular}
curve of degree $6$.

\

In dimension $3$ the situation is less satisfactory.
For Fano 3-folds we have the following general 

\begin{thm}
Rational points are potentially dense
on all Fano threefolds, with the possible exception of
the family of double covers of $\P^3$ 
ramified in a smooth sextic surface.
\end{thm}

This theorem follows from a case by case analysis
of Fano threefolds, classified by Fano, Iskovskikh,
Mori and Mukai (cf. \cite{mori-mukai-82}, \cite{mori-mukai-83}). 
Most threefolds on the list are unirational 
and rational points on unirational varieties are potentially
dense. The only exceptions for which unirationality 
is unknown for a general member of the corresponding family are:
smooth quartics in $\P^4$ (treated in \cite{harris-tschi-98}),
hypersurfaces of degree $6$ in the weighted projective space
$\P(1,1,1,2,3)$ (treated in this paper) and the double cover
of a smooth sextic in $\P^3$.

\

At the moment, there are no general approaches  
to the problem of Zariski--density of rational points, 
but a wealth of ad hoc methods appealing to specific
geometric features for different classes of varieties. 
In particular, the situation is such that 
one can prove that rational
points are potentially dense for 
certain  varieties which are higher in 
the Kodaira hierarchy, but not for all varieties which
are on the lower level. For example, density is well known
for abelian varieties, whereas it is an open problem
for general conic bundles over $\P^2$. 
It was pointed out to us by J.-L. Colliot-Th\'el\`ene that
Schinzel's hypothesis implies potential density for general 
conic bundles $X\ra \P^n$ (with $X$ smooth) (cf. 
Theorem 4.2, p. 73 in \cite{CT}).

\
 
A natural strategy is to look at fibrations with fibers abelian varieties
and to try to formulate sufficient conditions which would ensure 
that rational points are potentially dense. 
Our idea is roughly as follows: 
Let ${\cal A}\ra B$ be a semiabelian variety over
$B$  (defined over a number field $F$). 
Assume  further that there exists a 
subvariety ${\cal M}\subset {\cal A}$ which is a multisection
of the fibration and such that (potential) density holds for ${\cal M}$. 
Then there exists a finite extension of the groundfield $K/F$ such that
$K$-rational points are dense in a Zariski dense set of fibers
${\cal A}_b\subset {\cal A}$ (for $b\in B$). The rest of 
the paper is devoted to the detailed study of conditions 
and situations when this naive idea actually works.

\

{\small
{\bf Acknowledgments.} The second author is very grateful to 
Barry Mazur and Joe Harris for their interest, ideas and suggestions. 
We would like to thank J.-L. Colliot-Th\'el\`ene for his comments and 
the referee for the careful reading. 
The paper was written while both
authors were enjoying the hospitality of the Max-Planck-Institute in Bonn.
The first author was partially supported by the NSF.
The second author was partially supported by the NSA. 
}

\section{Elliptic fibrations}

Let $F$ be a field of characteristic zero and $B$ an irreducible $F$-scheme
of dimension $\ge 1$. Let $j \,:\,{\cal E}\ra B$ be a Jacobian elliptic fibration, 
i.e. an elliptic fibration with a zero section  $e\,:\, B\ra {\cal E}$, 
defined over $F$. 

\begin{defn}\label{dfn:non-torsion-section}
A {\em saliently ramified multisection} of ${\cal E}$ is 
an irreducible multisection ${\cal M}$ of 
$j$ which is defined over the algebraic closure 
$\overline{F}$ of $F$ and which is ramified in some non-singular fiber
${\cal E}_{b_0} = j^{-1}(b_0)$  $(b_0\in B(\overline{F}))$.
\end{defn}

\begin{prop} \label{prop:ell}
Let $j \,:\,{\cal E}\ra B$ be a Jacobian elliptic fibration, defined over
a number field $F$. 
If there exists a saliently ramified 
multisection ${\cal M}$ of ${\cal E}$
such that rational points on ${\cal M}$ are potentially dense
then rational points on ${\cal E}$ are potentially dense. 
\end{prop}

{\em Proof.}
Extending the groundfield $F$, if necessary, we can
assume that ${\cal M}$ is defined over $F$ and 
that $F$-rational points are dense on ${\cal M}$. 
Clearly, $j({\cal M}(F))$ is Zariski dense in $B$. 
We want to show that for almost all points $b\in j({\cal M}(F))$
(all outside a divisor in $B$) the points of intersection
$({\cal M}\cap {\cal E}_b)(F)$ generate an infinite
subgroup of the group of rational points ${\cal E}_b(F)$.

\begin{lm} \label{lemm:ell} 
Let $j\,:\,{\cal E}\ra B$ be a Jacobian elliptic fibration
satisfying the conditions of Proposition \ref{prop:ell}.
Then for any finite extension $K/F$ there exists a divisor
$D\subset B$ (depending on $K$), such that for all 
$b\in j({\cal M}(K))\cap (B\backslash D)(K) $ 
there exists a point
$p_b\in {\cal M}(K)\cap {\cal E}_b(K)$ which is
non-torsion in the group ${\cal E}_b(K)$.
\end{lm}

{\em Proof.} Extending the groundfield 
we can assume that the ramification  point $p_0$ of ${\cal M}$
in the  smooth fiber ${\cal E}_{j(p_0)}$ is $F$-rational  
and that ${\cal M}$  is defined over $F$.
For any $n$ the torsion points of order $n$
form a divisor $\Phi_n$ in ${\cal E}$. Thus, either our
irreducible multisection ${\cal M}$ is contained in  
$\Phi_n$ for some $n$, or it intersects each $\Phi_n$ in a proper subvariety, 
which projects onto a proper  subvariety  of the base $B$. 
Assume that the multisection ${\cal M}$ 
is contained in  $\Phi_n$  (for some $n$).
One can argue as follows:
fix an embedding of $\overline{F}\hookrightarrow \C$ and 
a neighborhood $U_0\subset B(\C)$ containing $b_0$
such that for all $b\in U_0$ the fibers ${\cal E}_b(\C)$ are smooth (this 
is possible by assumption). Then there exist a sequence of 
points $(b_i)_{i=1,...,\infty}$ converging to $b_0$ and
for each $b_i$ a pair of points $p_i,p_i'\in 
({\cal M}\cap {\cal E}_{b_i})(\C)$ (with $p_i\neq p_i'$) such that
both sequences $(p_i)_{i=1,..,\infty}$ and 
$ (p_i')_{i=1,..,\infty}$ converge to $p_0$. 
For each $i$ the difference $(p_i - p_i')\in {\cal E}_{b_i}(\C)$ is
torsion of order $n$ and 
at the same time it converges to $0=e(b_0)\in {\cal E}_{b_0}(\C)$.   
Contradiction. 

Consider the second possibility: 
Merel's theorem about the
boundedness of torsion for elliptic curves over number fields (cf. \cite{merel96})
implies that there exists a number $n_0(F)$ such that 
the  torsion subgroup of the group of rational 
points of every elliptic curve defined over $F$ is of order $\le n_0(F)$. 
Let $D$ be the union of the images of $\Phi_n\cap {\cal M}$ (over all $n\le n_0(F)$).
This is a divisor in $B$. 
Then for all $b\in j({\cal M}(F))$ which are not contained in 
$D$ there is an $F$-rational point $p_b$ in the fiber ${\cal E}_b$.
It cannot be torsion of order $\le n_0(F)$ on the one hand and 
there are no torsion points of order $>n_0(F)$. Therefore,
it must be a point of infinite order. 

\

An alternative argument would be to find a base change $B'\ra B$, \'etale at 
$j(p_0)\in B$ such that ${\cal M}$ pulls back to 
a {\em section} ${\cal M}'$ of ${\cal E}'\ra B'$. 
This section must be of infinite order in the Mordell-Weil group of ${\cal E}'$. 
Then one can apply specialization  (cf. \cite{silverman}).

\section{Double covers of $\P^2$}

In this section we work over a number field, allowing,
without further notice, finite extensions of the groundfield. 

\

Let $R\subset \P^2$ be a reduced curve of degree $6$ 
with a singular point $P$ and $\pi\,:\, D\ra \P^2$ a 
double cover of $\P^2$ ramified in $R$. 
Consider the blow up of $\P^2$ at $P$ 
together with a natural linear
system of lines in $\P^2$ passing through  $P$. 
We obtain an elliptic 
fibration ${\cal E}_{P}\ra \P^1$ on the 
desingularization ${\cal E}_P\ra D$ in $\pi^{-1}(P)$. 
The fibration ${\cal E}_P$ has
a natural section which is contained in the 
preimage of $P$ under this desingularization. 

\

\begin{prop}\label{prop:D}
Rational points on $D$ are potentially dense.
\end{prop}

{\em Proof.} Surely, this is so if ${\cal E}_{P}$ 
is a rational variety or 
birational to a $\P^1$-bundle over an elliptic curve.
This is the case if the curve $R$ contains multiple components or
if the point $P$ has  multiplicity $\ge 4$ on the curve
$R\subset \P^2$. 
Indeed, if $R$ contains multiple components, then
it is birational to a double cover of $\P^2$ ramified at
a curve of degree at most $4$. Then it is either rational
or birational to a $\P^1$-bundle over an elliptic curve.
Similarly, if the point $P$ has multiplicity $\ge 4$ at $R$ 
then the corresponding bundle ${\cal E}_{P}\ra \P^1$ is rational. 

\

In general,  ${\cal E}_{P}\ra \P^1$ is a
Jacobian elliptic fibration. We want to show the existence
of a saliently ramified elliptic or rational 
multisection in order to use the Proposition \ref{prop:ell}. 
We have to consider several possibilities, depending on 
the curve $R$. 

\begin{defn}(Tangent correspondence)
Let $R$ be a reduced curve in $\P^2$ and denote by $R^0$
the set of smooth points of $R$. 
The tangent correspondence ${\rm TC}(R)\subset R^0\times R^0$
is defined as follows: for any point $r\in R^0$ 
we consider the set of points ${\rm TI}_r(R)$
of transversal intersection of
the tangent line $\T_r(R)$  through $r$ with $R$. 
Then 
$$
{\rm TC}(R):= \{ (r, {\rm TI}_r(R))\}\subset R^0\times R^0.
$$ 
\end{defn}

\begin{lm} \label{lemm:tc}
Assume that ${\cal E}_{P}\ra \P^1$ is a
Jacobian elliptic fibration and that 
there exists a 
tangent correspondence ${\rm TC}(R)$ 
for the ramification curve $R$.
Assume that ${\rm TC}(R)$ 
surjects onto some open part of a component of $R^0$ 
(under the second projection). Assume in addition that 
we can choose such a component to be different from a line 
through $P\in R$. Then ${\cal E}_{P}\ra \P^1$ 
has a saliently ramified elliptic multisection. 
\end{lm}

{\em Proof.} The proof consists in a direct construction 
of two (different) elliptic curves 
contained in $D$  which are tangent
to each other at a point which is smooth in 
both curves and in $D$. One of these curves is going 
to be a fiber of the Jacobian fibration ${\cal E}_{P}\ra \P^1$
and the other a saliently ramified multisection of this fibration.

\

Denote by $L_{P,r}\subset \P^2$ the line
joining $P$ and $r$ (for $r\neq P$). 
We can find a pair of smooth points $r,r'\in R$ with
$(r,r')\in {\rm TC}(R)$ with the property:
The lines  $\T_{r}(R)$ and $L_{P,r'}$ both intersect in $r'$ and 
are both transversal to the tangent line $\T_{r'}(R)$
(by our assumptions on the tangent correspondence). We can 
also assume that the preimages of both lines are smooth
elliptic curves on $D$ (since the points $r'$ are allowed
to vary in a smooth family which is not contained in a line 
through $P$). These elliptic curves are tangent in $D$
at the unique preimage $\pi^{-1}(r')$. The proof of this
fact is local. Indeed, the image of the tangent
space  $\pi(\T_{\pi^{-1}(r')}(D))$ coincides with 
the tangent line $\T_{r'}(R)$.
Since the lines $\T_{r}(R)$ and $L_{P,r'}$  are transversal 
at $r'$ both elliptic curves $\pi^{-1}(\T_{r}(R))$
and $\pi^{-1}(L_{P,r'})$ are tangent at the point $\pi^{-1}(r')$ 
to the  $1$-dimensional subspace  which is the kernel 
of the differential ${\rm d}\pi$ at $\pi^{-1}(r')$.
Hence they are tangent at this point. 

\

The assumptions of \ref{lemm:tc} are satisfied
by most algebraic curves. More precisely

\begin{lm}\label{lemm:l}
Let $R\subset \P^2$ be any reduced curve with a point $P$ 
which contains either a component of degree $\ge 3$ 
or two conics or a conic and a line which doesn't contain $P$.
Then the assumptions of \ref{lemm:tc} are satisfied.
\end{lm}

{\em Proof.} Indeed, they are already satisfied for any
irreducible curve of degree $\ge 3$. A general tangent
line  at a smooth point $r\in R$ 
intersects the curve transversally
in at least one point $r'$ and the point $r'$ varies with $r$.
Thus, if a curve $R$ contains such components, we are done. 
The case of two conics and a conic and a line not containing
$P$ is evident.

\begin{lm}\label{lemm:2pts}
Let $L_1,...,L_6$ be a configuration of $6$ lines in $\P^2$
with no points of intersection of multiplicity $>3$. 
Then there are at least two points of intersection 
of multiplicity $2$. 
\end{lm}

{\em Proof.} Exercise.

\begin{cor} Rational points are potentially dense
on a double cover of any configuration $R$ of $6$ lines
in $\P^2$. 
\end{cor}

{\em Proof.} Consider two double points $P,Q$ and a line
$L\in R$ which doesn't contain $P$ and $Q$. 
Then we have two Jacobian elliptic fibrations ${\cal E}_P$ and 
${\cal E}_Q$ with fibers corresponding to the lines 
through $P$ and $Q$. Let $l\in L$ be a general
point. Then the lines $L_{P,l}$ and $L_{Q,l}$ are transversal
at $P$ to the line $L$. By the previous argument, the corresponding 
elliptic fibers of ${\cal E}_P$ and ${\cal E}_Q$ are tangent
at the preimage of $P$ in ${\rm d}\pi$. Hence, both Jacobian
elliptic fibrations have elliptic multisections tangent to  
smooth  fibers and we can use our standard argument \ref{prop:ell}.

\

This concludes the proof of \ref{prop:D}.

\section{Fano threefolds}

As described in the introduction, potential density 
of rational points on Fano threefolds
is proved by combining  classification and unirationality 
results. Smooth quartics for which 
unirationality is unknown were treated in \cite{harris-tschi-98}.

\

Let $V_1$ be the Fano threefold given as
a hypersurface of degree $6$ in the weighted projective
space $\P(1,1,1,2,3)$.
It can be realized as 
a double cover $V_1\ra V$ of the Veronese cone $V\subset \P^6$
ramified in an intersection of this cone with a
smooth cubic hypersurface $C\subset \P^6$ which doesn't
contain the vertex of $V$. We denote by $R$ the
ramification divisor.  If we blow up
the singular point of the cone $V$ we obtain a smooth
$3$-dimensional variety 
$V^* \simeq \P_{\P^2}({\cal O}(2)\oplus {\cal O}) $ which is birational to $V$.
This projective bundle over $\P^2 $ has two special sections
corresponding to the decomposition into the direct sum above.
We denote one of them by $\P_0:=\P_{\P^2}({\cal O}(2))$ and 
the other by $\P_{\infty}:= \P_{\P^2}({\cal O})$. The complement to
the singular vertex of $V$ is isomorphic to the line bundle
${\cal O}(2)$ over $\P_{\infty}$. The double cover of $V^*$ ramified
in $R$ (we will use the same notation $R$ for the 
preimage) will be denoted by $\pi\,:\, V_1^*\ra V^*$. We denote by 
$p_{\infty}$ the projection 
$p_{\infty} \,:\, R\ra \P_{\infty}$. On the reduced divisor
$R$ we can define a subvariety ${\rm disc}(R)$ by the condition
${\rm rk}({\rm d} p_{\infty}) <2$. This subvariety coincides with 
the subvariety of multiple intersection points of $R$ with
the fibers of the $\P^1$ bundle 
$\P_{\P^2}({\cal O}(2)\oplus {\cal O})$. 

\
 
We observe that $V_1^*$ has a structure of 
a Jacobian elliptic fibration
over $V_1^*\ra \P^2_{\infty}$ where the elliptic fibers project
to the lines of the cone $V^*$ and the zero section is
$\P_0=\P_{\P^2}({\cal O}(2))$. In this picture the preimage 
of ${\rm disc}(R)$ coincides with singular points of 
the fibers of the elliptic fibration $V_1^*\ra \P^2_{\infty}$.

\begin{thm}
Rational points on $V_1$ are potentially dense. 
\end{thm}

{\em Proof.}
The Jacobian elliptic fibration 
$V_1^*\ra \P^2_{\infty}$  has many
double sections which are K3 surfaces.
Any non-zero section $s$ of ${\cal O}(2)$ over $\P_{\infty}$
determines an embedding of 
${\cal O}\hookrightarrow  {\cal O} \oplus {\cal O}(2)$ and 
therefore a section $\P^2_s\simeq \P^2\hookrightarrow V^*$. 
The section $\P^2_s$ intersects the 
ramification divisor $R\subset V^*$ 
along a curve $R_s\subset \P^2_s$ of degree $6$. 
Thus we obtain a surface $D_s\subset V_1^*$ 
which is a double cover of $\P^2_s$ ramified in $R_s$. 
We will show that among the $D_s$ we can find one
which is a saliently ramified multisection of the Jacobian 
elliptic fibration $V_1^*\ra \P^2_{\infty}$  and 
on which rational points are potentially dense. 
Then we apply \ref{prop:ell}.

\

We choose a smooth point $P$ on the reduced 
ramification divisor $R$ 
such that the differential
of the projection to the base $\P^2$ of 
$\P_{\P^2}({\cal O} \oplus {\cal O}(2))$ is an isomorphism. 
This is always possible unless $R$ is a $\P^1$-fibration
over a curve of degree $6$ in $\P^2$. But in this case
the initial cubic hypersurface $C\subset \P^6$ must contain
the vertex of the Veronese cone $V$. It is worth noticing that
the corresponding smooth threefold is a $\P^1$-fibration over
a K3 surface and therefore not rationally connected. 
By the same reasoning, ${\rm disc}(R)$ is a curve in $R$. 

\

Continuing the proof, we can now guarantee the 
existence of a family of sections
${\cal F}_P \simeq H^0(\P^1, {\cal O}(2))$
such that for all $s\in {\cal F}_P$ the
section $\P^2_s$ is tangent to the 
ramification divisor $R\subset V^*$ at $P$. 
This is possible because we can find 
a quadratic function on $\F^2$ with given 
value and given first differential at $P$, which defines
the desired section $s$ of ${\cal O}(2)$ on $\P_{\infty}$.
Moreover, locally around $P$
the second differential of the function describing $R$ with
respect to $\P^2_s$ can be chosen at will.  
Then the curve $R_s$  is a singular curve of 
degree $6$ with singular point $P$. This implies
that the corresponding double cover $D_s$ is a singular
surface with an isolated singular point at the preimage of $P$.
Moreover, notice that $P\not\in {\rm disc}(R)$.
Therefore, in every curve $R_s$ 
there exists an open neighborhood $U_s(P)\subset R_s$ of $P$ 
which doesn't intersect the curve ${\rm disc}(R)\subset R$.
This shows that every such $D_s$ is a saliently ramified multisection
of the fibration $V_1^*\ra \P^2_{\infty}$.

\

Now we want to use this multisection which is tangent to
points contained in smooth fibers of the Jacobian elliptic fibration 
$V_1^*\ra \P^2_{\infty}$ along $U_s(P)$ in our argument \ref{prop:ell}. 
We need to show that rational points are
potentially dense on ${\cal E}_{P,s}$ (or, equivalently on $D_s$).
But this is the content of Proposition \ref{prop:D}. 
Therefore, rational points are potentially dense on $V_1$.

\begin{rem} In fact, we have a choice of $D_s$. In particular,
{\em in most} situations we can choose $D_s$ such that the ramification
curve $R_s$ has two singular points of multipliticity two. 
Then we obtain two Jacobian elliptic fibrations with some non-singular
fibers of one fibration tangent to non-singular fibers of the other.
\end{rem}

\nocite{mazur-92}
\nocite{manin-93}

\end{document}